\documentclass{article}
\usepackage{amsmath,amsthm,amsfonts,latexsym}

\begin{document}

\newtheorem*{theo}{Theorem}
\newtheorem*{pro} {Proposition}
\newtheorem*{cor} {Corollary}
\newtheorem*{lem} {Lemma}
\newtheorem{theorem}{Theorem}[section]
\newtheorem{corollary}[theorem]{Corollary}
\newtheorem{lemma}[theorem]{Lemma}
\newtheorem{proposition}[theorem]{Proposition}
\newtheorem{conjecture}[theorem]{Conjecture}
%{\theorembodyfont{\rmfamily}
 \newtheorem{definition}[theorem]{Definition}
 \newtheorem{remark}[theorem]{Remark}
%}
\newcommand{\Naturali}{{\mathbb{N}}}
\newcommand{\Reali}{{\mathbb{R}}}
\newcommand{\Complessi}{{\mathbb{C}}}
\newcommand{\Toro}{{\mathbb{T}}}
\newcommand{\Relativi}{{\mathbb{Z}}}
\newcommand{\HH}{\mathfrak H}
\newcommand{\KK}{\mathfrak K}
\newcommand{\LL}{\mathfrak L}
\def\A{{\cal A}}
\def\B{{\cal B}}
\def\H{{\cal H}}
\def\L{{\cal L}}
\def\N{{\cal N}}
\def\M{{\cal M}}
\def\O{{\cal O}}
\def\U{{\cal U}}

\title{
On infinite tensor products of projective \\
unitary representations}

\author{Erik B\'edos$^*$,
Roberto Conti$^{**}$ \\}
\date{December 11, 2000}
\maketitle
\markboth{R. Conti, Erik B\'edos}{Infinite tensor products and projective
representations}
\renewcommand{\sectionmark}[1]{}
\begin{abstract}
We initiate
a study of infinite tensor products of
projective unitary representations of a discrete group $G$.
Special attention is given
to regular representations twisted by $2$-cocycles
and to projective representations associated with CCR-representations of
bilinear maps. Detailed computations
 are presented in the case where $G$ is a finitely generated free
abelian group.
We also  discuss an extension problem about product
type actions of $G$,
where
the projective representation theory of $G$
plays a central role.

\vskip 0.3cm
\noindent {\bf MSC 1991}: 22D10, 22D25, 46L55, 43A07, 43A65

\smallskip
\noindent {\bf Keywords}: projective unitary representation, infinite tensor
product,
amenable group.

\bigskip
\end{abstract}

\vfill
\thanks{\noindent $^*$  partially supported by the Norwegian Research
Council.\par
 \noindent Address: Institute of Mathematics, University of
Oslo,
P.B. 1053 Blindern,\\ 0316 Oslo, Norway. E-mail: bedos@math.uio.no. \\

\noindent $^{**}$ supported by the EU TMR network ``Non Commutative
Geometry.'' \\
Adress: Dipartimento di Matematica, Universit\`a di Roma
``Tor Vergata'',\\
I-00133 Roma, Italy. E-mail: conti@mat.uniroma2.it.\par}

\newpage

\section{Introduction}
The theory of infinite tensor products of Hilbert spaces started
with the seminal paper by von Neumann \cite{vN}.
Later on, Guichardet \cite{Gui,Gui2} approached this matter from a
slightly different point of view and
developed a unified framework for treating
several related  concepts involving operators, algebras and
functionals.
The notion of infinite tensor product has been
mainly used in this form in operator algebras and quantum field theory
over the last three decades (see \\ e.\ g.\ \cite{EK} for a recent
overview).

The existence of some infinite tensor product of representations of a group
has been established and used in some recent works. For example, it
was shown in \cite{ABC} that a locally compact group is
$\sigma$-compact and amenable
if and only if  there exists
an infinite tensor power of its regular representation. Such an
infinite tensor power construction was then a useful tool for studying
covariance of certain (induced) product-type representations of
generalized Cuntz algebras
with respect to natural product-type actions.
This circle of ideas has been generalized and thoroughly investigated in
\cite{BC}.
In another direction, the infinite tensor product of certain unitary
representations of some group of diffeomorphisms was shown to exist
under suitable assumptions
in \cite{HiSh}.

In this paper we initiate a study of
infinite tensor products of \emph{projective} unitary representations of a
 discrete group. It is in fact not obvious that such infinite tensor
products
 exist at all. Indeed it is quite easy to realize that it is
impossible to form the infinite tensor power of a single projective unitary
representation unless the associated 2-cocycle vanishes.  Besides its
intrinsic interest, this
new generality has
the potential advantage to allow for extensions of the analysis given in
\cite{ABC,BC} to a
broader class of product-type actions on
the $0^{th}$-degree part of extended Cuntz algebras.
It is also relevant when studying extensions of
product-type actions from the algebraic to the von Neumann algebra level.
Finally it provides a way to represent faithfully on
infinite tensor product spaces some familiar C*-algebras like
non-commutative tori.
To avoid technicalities, we stick to the case of a discrete
group, although it could be of interest  in the future to consider a
locally compact (or even just a topological)
group and strongly continuous projective unitary representations
of such a group.

The paper is organized as follows.
Section 2 is devoted to some preliminaries on projective unitary
representations, product sequences of 2-cocycles
and infinite tensor products.
Section 3 contains our main existence results
for infinite tensor
products of projective unitary representations. We especially display some
sufficient conditions for countable amenable groups in the case of
projective
regular representations
and in the case of projective representations associated with
CCR-representations
 of bilinear maps.
To illustrate our work
we specialize
in section 4 to the case
of finitely generated free abelian groups.
The final section deals with infinite tensor
products of actions of a discrete group $G$ on von Neumann algebras.
We concentrate our attention to the existence problem of such product
actions in the case of unitarily implemented
actions.
One of our result exhibits an obstruction
for extending some algebraic tensor power action of $G$ to the weak
closure that lies in the second cohomology group $H^{2}(G,\Toro)$. In
another result, the obstruction lies in the non-amenability of $G$.

\section{Preliminaries}
Throughout this note
$G$ denotes a  non-trivial discrete group, with neutral element $e$.

A \emph{2-cocycle} (or multiplier) on $G$ with values in the circle group
$\Toro$ is a map
$u: G \times G \to \Toro$
such that
$$u(x,y)u(xy,z)=u(y,z)u(x,yz) \quad\quad (x,y,z \in G) ,$$
see e.g. \cite[Chapter IV]{Br}.
We will consider only \emph{normalized} 2-cocycles, satisfying
$$u(x,e)=u(e,x)=1 \quad\quad (x \in G).$$
The set of all such 2-cocycles, which is denoted by $Z^2(G,\Toro)$,
becomes an abelian group under pointwise product. We  equip
$Z^2(G,\Toro)$ with the topology of pointwise convergence.

A 2-cocycle $v$ on $G$ is called a \emph{coboundary} whenever
$v(x,y)=\rho(x)\rho(y)\overline{\rho(xy)}$ $(x,y \in G)$
for some $\rho: G \to \Toro$, $\rho(e)=1$, in which case we write
$v = {\rm d}\rho$
(such a $\rho$ is uniquely determined up to multiplication by a
character). The set of all coboundaries, which is denoted by
$B^2(G,\Toro)$, is a  subgroup of $Z^2(G,\Toro)$, which is easily
seen to be closed.
(Indeed, assume that $({\rm d}\rho_{\alpha})$ is a net in
$B^2(G,\Toro)$ converging to $v \in Z^2(G,\Toro).$ Due to Tychonov's
theorem, we may, by passing to a subnet if necessary, assume that
$\rho_{\alpha}$ converges pointwise to $\rho$, for some
$\rho: G \to \Toro$, $\rho(e)=1.$ Then we have $v = {\rm d}\rho.$)

The quotient group $H^2(G,\Toro):=Z^2(G,\Toro)/B^2(G,\Toro)$
is called the \emph{second cohomology group} of $G$ with values
in $\Toro$. We denote elements in $H^2(G,\Toro)$ by $[u]$ and
write $v \sim u$ when $[v] = [u] \quad (u, v \in Z^2(G,\Toro))$. We also
write $v \sim_{\rho} u$  when  we have $v = ({\rm d}\rho)  u$
for some coboundary ${\rm d}\rho$.

\vspace{1ex}
We recall a few facts concerning infinite products of complex
numbers (see \cite{vN}).
Let $(z_{i})$ denote a sequence of complex numbers. We say that the
infinite product $\prod_{i} z_{i}$ exists (or converges) if the limit
of the net  $(\prod_{i \in J} z_{i}) _{J \in  \cal{F}}$ exists,
where $\cal{F}$ denotes the family of non-empty finite subsets of
$\Naturali$ ordered by inclusion. We then also use $\prod_{i} z_{i}$
to denote this limit. We will need the following result:

\vspace{1ex}
\noindent Assume that $\sum_{i} | 1 - z_{i} | < \infty.$ Then $\prod_{i}
z_{i}$
exists, and $\prod_{i} z_{i} \neq 0$ if all $z_{i}$'s are non-zero.
Conversely, assume that $\prod_{i} z_{i}$ converges to a non-zero
element. Then $\sum_{i} | 1 - z_{i} | < \infty.$

\vspace{1ex}
We shall be  interested in \emph{product} sequences in $Z^2(G,\Toro)$:
we call a
sequence $(u_{i})$ in $Z^2(G,\Toro)$
a \emph{product} sequence whenever the (pointwise) infinite product
$ u = \prod_{i} u_{i}$ exists on $ G \times G$
($u$ being then
obviously a 2-cocycle itself).

A cohomological problem concerning product
sequences  is that  perturbing a product sequence (by
a coboundary in each component) does not necessarily lead to a
product sequence, as may be illustrated by taking all $u_{i}$'s to be
1 and perturbing by the same coboundary $v \neq 1$ in each component.
The following lemma  somewhat clarifies this problem.

\begin{lemma} Let $(u_{i})$ and $(v_{i})$ be two sequences in
$Z^2(G,\Toro)$ satisfying $ v_{i}\sim_{\rho_{i}} u_{i}$ for every $i$.

\vspace{1ex}
\noindent i) Assume that $\rho := \prod_{i} \rho_{i}$ exists. Then $(v_{i})$
is
a product sequence if and only if $(u_{i})$ is a product sequence, in
which case we have $\prod_{i} v_{i} \sim_{\rho} \prod_{i} u_{i}$.

\vspace{1ex}
\noindent
ii) Assume that $(u_{i})$ and $(v_{i})$ are both product sequences.
Then $\prod_{i} v_{i} \sim \prod_{i} u_{i}$ (even if $\prod_{i}
\rho_{i}$ does not necessarily exist).
\end{lemma}

\begin{proof}
As i) is straightforward, we only show ii). So we assume that
$u = \prod_{i} u_{i}$ and $v = \prod_{i} v_{i}$ both exist. Then $
w:= \prod_{i} {\rm d}\rho_{i} =  \prod_{i} \overline{u_{i}} v_{i}$ also
exists and is the limit of a net of 2-coboundaries. As
$B^2(G,\Toro)$ is closed, this implies that
$w \in B^2(G,\Toro)$. Since $ v = w u$, this shows that $v \sim u$, as
asserted.

\noindent  (To see that $\prod_{i}
\rho_{i}$ does not necessarily exist, assume that $G$ possess a non-trivial
character $\gamma.$ Set $u_{i}=v_{i}=1$ and $\rho_{i} = \gamma$ for
all $i.$ Then
clearly $ v_{i}\sim_{\rho_{i}} u_{i}$ while $\prod_{i}
\rho_{i}$ does not exist.)
\end{proof}

\vspace{1ex}

A \emph{projective unitary representation} $U$ of $G$ on a Hilbert space
$\mathcal{H}$ associated with some $u \in Z^2(G,\Toro)$ is a map from $G$
into the group of unitaries on $\mathcal{H}$ such that
$$ U(x)U(y) = u(x,y) U(xy) \quad\quad (x,y \in G).$$
If we pick a $\rho: G \to \Toro$ satisfying $\rho(e)=1$
and set $V = \rho \, U$, then $V$ is also a projective unitary
representation of $G$ on $\mathcal{H}$
associated with a 2-cocycle $v$ satisfying $v \sim_{\rho} u$.
Such a $V$ is called a \emph{perturbation} of $U$.

\vspace{1ex}
To each  $u \in Z^2(G,\Toro)$ one may associate the  left $u$-regular
projective
unitary representation $\lambda_{u}$ of $G$ on $\ell^2(G)$  defined  by
$$(\lambda_{u}(x)f)(y)=u(y^{-1},x)f(x^{-1}y) \quad\quad (f \in
\ell^2(G),\ x,y \in G).$$
Choosing $u=1$  gives  the left regular representation of $G$
which we will just denote by $\lambda$. It is well known (and easy to see)
that if $v \sim_{\rho} u$, then
$\lambda_{v}$ is unitarily  equivalent to $ \rho  \lambda_{u}$.

\vspace{1ex}
For $i = 1, 2$, let $U_{i}$ be a projective unitary representation of $G$
on a Hilbert space $\mathcal{H}_{i}$
associated with   $u_{i} \in Z^2(G,\Toro)$. Then the naturally
defined tensor product
representation $U_{1} \otimes U_{2}$ is easily seen to be a projective
unitary representation of $G$ on the Hilbert space
$\mathcal{H}_{1} \otimes \mathcal{H}_{2}$
associated with the product cocycle  $u_{1} u_{2}$. In the case of
ordinary unitary representations of a group, it is a classical result
of Fell (cf. \cite{Di}, 13.11.3) that the left regular representation
acts in an absorbing way with respect to tensoring (up to multiplicity and
equivalence). In the projective case we have the following analogue.

\medskip
\begin{proposition}
Let $u, v$ be elements in $Z^2(G,\Toro)$ and let $V$ be any projective
unitary
representation of $G$ on a Hilbert space $\mathcal{H}$
associated with $v$. Then the tensor product representation $\lambda_{u}
\otimes V$ is unitarily equivalent to $\lambda_{uv} \otimes
id_{\mathcal{H}}$, i.e.
to (dim $V) \cdot \lambda_{uv}$.
\end{proposition}

\begin{proof} We leave to the reader to check that the
same unitary operator $W$ as in the non-projective case ( which is
determined on $\ell^2(G) \otimes \mathcal{H} \ (\, \cong
\ell^2(G,\mathcal{H}))$ by
$(W (f \otimes \psi)) (x) = f(x) V(x^{-1}) \psi$) implements the asserted
equivalence.
\end{proof}

\medskip

We conclude this  section with a short review on infinite tensor products of
Hilbert spaces
and operators.
(See \cite{Gui,Gui2} for more information.)

Let $\mathcal{H} = \{\mathcal{H}_{i}\}$ denote  a sequence of Hilbert
spaces and  $\phi = \{\phi_{i}\}$ be  a sequence of unit vectors
where $\phi_{i} \in \mathcal{H}_{i}$ for each $i \geq 1$. We denote by
$\mathcal{H}^{\phi}$ or  by $\bigotimes_{i}^{\phi} \mathcal{H}_{i}$
the associated infinite tensor product
Hilbert space of the $\mathcal{H}_{i}$'s along the sequence $\phi$.

 For any sequence $\psi_{i} \in \mathcal{H}_{i}$ such that
$$\sum_i |\,1- \| \psi_i \| \,| < \infty \; \text{and} \;  \sum_i
|\,1-(\psi_i,\phi_{i})\,| < \infty, $$
there corresponds a so-called decomposable vector in $\mathcal{H}^{\phi}$
denoted by $\otimes_{i} \psi_{i}$.
 If $\otimes_{i} \xi_{i}$ is
another decomposable vector in $\mathcal{H}^{\phi}$, then
$$ (\otimes_{i} \psi_{i} , \otimes_{i} \xi_{i}) =
\prod_{i} (\psi_{i} ,  \xi_{i}) .$$
A decomposable vector of the form $\psi_{1} \otimes \ldots \otimes
\psi_{k} \otimes \phi_{k+1} \otimes \phi_{k+2} \otimes \ldots$
is called elementary. The set of all elementary decomposable vectors
is total in  $\mathcal{H}^{\phi}$.

Let $T_{1},T_{2},\ldots$ be  a sequence of bounded linear operators where
each $T_{i}$
acts on $\mathcal{H}_{i}$. For each fixed $n \in \Naturali$ there
exists a unique bounded linear operator $\tilde{T}_{n}$ acting on
$\mathcal{H}^{\phi}$ which is determined by
$$ \tilde{T}_{n} (\otimes_{i} \psi_{i}) = T_{1}\psi_{1} \otimes \ldots
\otimes T_{n}\psi_{n} \otimes \psi_{n+1}\otimes \psi_{n+2} \otimes
\ldots $$
for each decomposable vector $\otimes_{i} \psi_{i}$.
Similarly, one may define $\tilde{T}_{J}$ for each (nonempty) finite
$J \subset  \Naturali.$
 Under certain
assumptions, the net $\{\tilde{T}_{J}\}$ converges in the strong
operator topology to a bounded linear operator on
$\mathcal{H}^{\phi}$ which is then denoted by $\otimes_{i} T_{i}$.

By \cite[Part II, Proposition 6]{Gui2}), a sufficient condition for
$\otimes_{i} T_{i}$ to exist  is that
$$ \prod_{i} \| T_{i} \|  \ \mbox{exists} \, , \,
\sum_{i} |1 - \|T_{i} \phi_{i}\| \, | < \infty  \ \mbox{and} \
\sum_{i} |1 - (T_{i} \phi_{i}, \phi_{i}) \, | < \infty, $$
in which case we
have $(\otimes_{i} T_{i}) \, (\otimes_{i} \psi_{i}) \, = \, \otimes_{i}
T_{i}\psi_{i}$ for all elementary decomposable vectors $\otimes_{i}
\psi_{i}.$\\

When all $T_{i}$'s are unitaries (which is the case of interest
in this paper) we have the following
result, which will be used several times in the sequel.

\begin{proposition} Let $(T_{i})$ be a sequence of unitaries where
each $T_{i}$
acts on $\mathcal{H}_{i}$. Then $\otimes_{i} T_{i}$ exists on
$\mathcal{H}^{\phi}$ if and only if
$$ (*) \ \ \ \sum_i |1-(T_{i}\phi_i,\phi_{i})| < \infty,$$
in which case $\otimes_{i} T_{i}$ is a unitary on $\mathcal{H}^{\phi}$
satisfying $(\otimes_{i} T_{i})^{*} = \otimes_{i} T_{i}^{*}.$
\end{proposition}

\begin{proof} Assume first that $(*)$ holds. It is then quite elementary to
deduce from Guichardet's result mentioned above  that $\otimes_{i}
T_{i}$ and $\otimes_{i} T_{i}^{*}$ both exist.
 Moreover, these two operators are then  isometries, being
 the strong limit of a net of unitaries, and they are easily seen to be the
inverse
of each other.
 So both are  unitaries satisfying $(\otimes_{i} T_{i})^{*} = \otimes_{i}
T_{i}^{*}.$

 Assume now that $T:= \otimes_{i} T_{i}$ exists on
 $\mathcal{H}^{\phi}.$ Then $T$ is non-zero (being an isometry), so
 there are elementary decomposable vectors $\otimes_{i} \psi_{i}$ and
 $\otimes_{i} \xi_{i}$
 such that
$$0 \neq c:=(T \, \otimes_{i} \psi_{i} , \otimes_{i} \xi_{i}).$$
Let $J$ be
any finite subset of $\Naturali$ large enough so that $\psi_{i} = \xi_{i} =
\phi_{i}$ for all $i \notin J.$ Then we have
$$(\tilde{T}_{J} \otimes_{i} \psi_i , \otimes_{i} \xi_i) =
\prod_{i \in J} (T_{i} \,  \psi_{i} ,  \xi_{i}).$$
Since $T=\lim_J \tilde{T}_J$, we get
$c=\lim_J \prod_{i \in J} (T_{i} \,  \psi_{i} ,  \xi_{i})$,
 i.e. $\prod_{i \in \Naturali} (T_{i} \,  \psi_{i} ,  \xi_{i})$
converges
to a non-zero value.

\noindent Thus we get $\sum_i | 1-(T_{i} \,  \psi_{i} ,  \xi_{i}) | <
\infty$
and therefore
$\sum_i | 1-(T_{i} \,  \phi_{i} ,  \phi_{i}) | < \infty $ since
$\psi_{i} = \xi_{i} = \phi_{i}$ for all but finitely many $i$'s.
\end{proof}

\section{Infinite tensor products of projective unitary
representations}

 Before attacking the main problem whether it is possible to form an
infinite
tensor product of  a sequence of projective unitary representations,
at least in some cases, we first show that this
 construction, when possible, produces a new projective unitary
 representation of $G$, and also make some general observations.

 \begin{theorem} Let $U_{i}$ be a sequence of projective unitary
 representations of G acting respectively on a Hilbert space
 $\mathcal{H}_{i}$ and with associated $u_{i} \in Z^2(G,\Toro)$.
 Let $\phi = (\phi_{i})$ be a sequence of unit vectors where each
 $\phi_{i} \in \mathcal{H}_{i}$. Assume that $\otimes_{i} U_{i}(x)$
 exists on $\mathcal{H}^{\phi} = \otimes_{i}^{\phi}\mathcal{H}_{i}$
 for each $x \in G$. Then we have

 \vspace{1ex} \noindent
i)  $(u_{i})$ is a product sequence in $Z^2(G,\Toro).$

 \vspace{1ex} \noindent
ii) The map $x \rightarrow U^{\phi}(x):= \otimes_{i} U_{i}(x)$ is a
 projective unitary representation of G on $\mathcal{H}^{\phi}$ with
  $u = \prod_{i} u_{i}$ as its associated 2-cocycle.

 \vspace{1ex} \noindent
iii) If there exists one $k$ such that $U_{k}$ is unitarily
  equivalent to $\lambda_{u_{k}}$, then
$U^{\phi}$ is unitarily equivalent to $\lambda_{u} \otimes
id_{\mathcal{H}}$, where $\mathcal{H}$ denotes any infinite
dimensional separable Hilbert space.

 \vspace{1ex} \noindent
iv) $\lambda \otimes U^{\phi}$  is unitarily equivalent to $\lambda_{u}
\otimes id_{\mathcal{H}^{\phi}}$.

 \end{theorem}

 \begin{proof}Notice first that  Proposition 2.3  implies that
each $
 U^{\phi}(x):= \otimes_{i} U_{i}(x)$ is a unitary.

 i) Let $g,h \in G$. We must show that $\prod_{i} u_{i}(g,h)$ converges.
Now

 $$\otimes_{i} U_{i}(gh)$$ and
 $$(\otimes_{i} U_{i}(g)) (\otimes_{i}
 U_{i}(h)) = \otimes_{i} U_{i}(g)U_{i}(h) = \otimes_{i}
 u_{i}(g,h)\,U_{i}(gh)$$  are both unitaries.
Putting $a_{i} = (\, U_{i}(gh))\phi_{i} \, ,\, \phi_{i} \,)$, we deduce
from Proposition 2.3  that
 $$\sum_{i} |\,1 - a_{i} \,| < \infty \quad \text{and} \quad
 \sum_{i} |\,1 - u_{i}(g,h) a_{i} \,| < \infty .$$
This implies that
$\sum_{i} |\,1 - u_{i}(g,h) \,| < \infty $, and therefore that
$\prod_{i} u_{i}(g,h)$ converges, as desired. ( We use here
implicitely that
whenever $z \in \Toro$ and $a \in \Complessi$, then $|1-z| =
|1-\overline{z}| \leq |1-a| + |a-\overline{z}| = |1-a| + |za-1|$).

ii) Using i) we get
$$ U^{\phi}(x) \, U^{\phi}(y) = \otimes_{i}
 u_{i}(x,y)\,U_{i}(xy) = (\prod_{i} u_{i}(x,y)) \, \otimes_{i} U_{i}(xy) =
 u(x,y)\, U^{\phi}(xy)$$
 for all $x,y \in G$, as asserted.

 iii) and iv) follow easily from Proposition 2.2.
\end{proof}

An obvious, but noteworthy consequence of part i) of this theorem is that it
is
impossible to form the infinite tensor power of a single projective unitary
representation unless the associated 2-cocycle vanishes.
In another direction, the case where infinitely many of the $U_{i}$'s are
projective  regular
representations of $G$ can not occur in this theorem when $G$ is uncountable
or
non-amenable, as easily follows from our next theorem.
(We refer to \cite{Pat}  or \cite{Pi} for information on
amenability).

\begin{theorem} Let $(u_{i})$ be a sequence in $ Z^2(G,\Toro)$ and set
 $U_{i} = \lambda_{u_{i}}$ for every $i$.
 Let $\phi = (\phi_{i})$ be a sequence of unit vectors in $\ell^2(G).$
 Assume that $\otimes_{i} U_{i}(x)$
 exists on $\mathcal{H}^{\phi} = \otimes_{i}^{\phi}\ell^2(G)$
 for each $x \in G$. Then $G$ is countable and amenable.
 \end{theorem}

 \begin{proof}
Using Proposition 2.3, it follows  that
$\sum_i |1-(U_{i}(x)\phi_i,\phi_i)|<\infty$ for
every $x \in G$.
Notice that
$$|(U_{i}(x)\phi_i,\phi_i)|
\leq
(\lambda(x)|\phi_i|,|\phi_i|)
\leq 1.$$
Hence we get
$$(\lambda(x)|\phi_i|\, , \,|\phi_i|) \to 1 \quad\quad (x \in G).$$
This means that  the trivial 1-dimensional representation of $G$ is weakly
contained in $\lambda $ and the amenability of $G$ follows.

\noindent Moreover, setting
$f_i(x):=|(\lambda(x)\phi_i,\phi_i)| \geq 0$
we have $0 \leq f_i \leq 1$,
$f_i \in C_0(G)$ (cf. \cite[13.4.11]{Di})  and $f_i \to 1$ pointwise.
Then $f_i^{-1}([1/2,1])=:H_i$ is finite and $G = \cup_i H_i$, so $G$ is
countable.
\end{proof}

\vspace{2ex}
In view of this theorem, it is quite natural to wonder whether some
converse holds. We shall provide a partial answer in Corollary 3.4.
To ease our exposition, we introduce some terminology.
A sequence $(F_{i})$
of non-empty, finite subsets of $G$ will be called a
$F$-sequence (resp.\ $\sigma F$-sequence) for $G$ whenever
$$\lim_{i} \frac{\# (F_{i} \cap xF_{i})}{\# F_{i}} = 1 \ \text{for
all}\  x \in G,$$
$$\text{(resp.} \sum_{i}|1 - \frac{\# (F_{i} \cap xF_{i})}{\# F_{i}}| <
\infty
\ \text{for all} \ x \  \in G).$$
A $F$-sequence $(F_{i})$ for $G$ is often called a \emph{F\"{o}lner}
sequence
in the literature. We remark that the definition is usually phrased
in a slightly different, but equivalent, way (involving the symmetric
difference of sets) and that some authors also require that $F_{i}
\subseteq F_{i+1}$ for every $i$. Anyhow, thanks to F\"{o}lner (see
\cite{Pat,Pi}), we know
that  $G$ is countable and amenable if and onl y if $G$ has a $F$-sequence.
Now, obviously, a  $\sigma F$-sequence for $G$ is also a $F$-sequence.
Moreover,
any $F$-sequence has some subsequence which is a $\sigma F$-sequence,
as is easily checked. Hence we can also conclude that  $G$ is
countable and amenable  if and only if $G$ has a $\sigma F$-sequence.

\vspace{2ex}
When  $F$ is a  subset of $G$, we denote by $\chi_F$ its
characteristic  function.

\begin{theorem} Let $(u_{i})$ be a sequence in $Z^2(G,\Toro).$ Assume that
$G$ is countable and amenable, and has a $\sigma F$-sequence $(F_{i})$ which
satisfies
$$ (*)  \ \ \sum_{i} \frac{1}{\# F_{i}} \sum_{y \in F_{i}}|1 -
u_{i}(y^{-1},x) | < \infty
 \ \text{for all} \ x \in G.$$
Set $U_{i} = \lambda_{u_{i}}$ and $\phi_i:=\chi_{F_{i}}/{(\#
 F_{i})}^{1/2}$ for every $i$.\\
 Then $\phi = (\phi_{i})$ is a sequence
 of unit vectors in  $\ell^2(G)$ such that $\otimes_{i} U_{i}$
 exists on $\mathcal{H}^{\phi} = \otimes_{i}^{\phi}\ell^2(G)$.
 \end{theorem}

 \begin{proof} We first record  some easy calculations. Let  $F$ be a finite
 (non-empty) subset  of  $G$ and set  $\phi_F:=\chi_{F}/{(\#
 F)}^{1/2}$. Let $u \in Z^2(G,\Toro)$.
  Then we have
  $$ (\lambda(x)\phi_F,\phi_F)=\frac{1}{\# F}\# (F \cap xF)$$
 for every $x \in G.$ More generally we have

  $$(\lambda_u(x)\phi_F,\phi_F)
=\frac{1}{\# F}\sum_{y \in F \cap xF}u(y^{-1},x) \ $$

and therefore
$$(\,(\lambda(x) - \lambda_u(x))\phi_F,\phi_F\,) =
\frac{1}{\# F}\sum_{y \in F \cap xF}(1 - u(y^{-1},x)) \ $$
for all $x \in G$. \\

\vspace{2ex}
Using the triangle inequality and the above computations, we get

$$ \ \sum_i |1-(U_{i}(x)\phi_i,\phi_i)| \leq \sum_i
   |1-( \lambda(x)\phi_i,\phi_i)| + \sum_i
    |(\, (\lambda(x) - U_{i}(x)) \phi_i,\phi_i \,)|$$

$$=\sum_{i}|1 - \frac{\# (F_{i} \cap xF_{i})}{\# F_{i}}|  \  + \ \sum_{i}
\frac{1}{\# F_{i}} \ |\sum_{y \in F_{i} \cap
xF_{i}}(1 - u_{i}(y^{-1},x) ) \ | $$
$$ \leq \sum_{i}|1 - \frac{\# (F_{i} \cap xF_{i})}{\# F_{i}}|  \  + \
\sum_{i}
\frac{1}{\# F_{i}} \ \sum_{y \in F_{i}}|1 -
u_{i}(y^{-1},x) | $$
 for all $x \in G.$ Since  $(F_{i})$ is a $\sigma F$-sequence for
 $G$ satisfying $(*)$, both sums above converge for all $x\in G$.
 Hence, $ \ \sum_i |1-(U_{i}(x)\phi_i,\phi_i)| < \infty $ for all $x\in G$
 and the assertion follows from Proposition 2.3.
 \end{proof}

Clearly, if $u_{i}=1$ for all but finitely many $i$'s, any $\sigma
F$-sequence
$(F_{i})$ for $G$ trivially satisfies $(*)$. In this
case, the above theorem could also have been deduced from \cite{C2}.

 \begin{corollary} Let $G$ be countable and amenable, and let $(v_{j})$ be a
 product sequence in $Z^2(G,\Toro)$. Then there exist
 a subsequence $(u_{i})$ of $(v_{j})$ and a sequence $\phi =
 (\phi_{i})$ of unit vectors in  $\ell^2(G)$ such that $\otimes_{i}
\lambda_{u_{i}}$
 exists on $\mathcal{H}^{\phi} = \otimes_{i}^{\phi}\ell^2(G)$.
 \end{corollary}

  \begin{proof} First we pick  a $\sigma F$-sequence $(F_{i})$ for
  $G$ and a growing sequence $(H_{i})$ of non-empty finite subsets
  of $G$ satisfying $\cup_i H_i = G$. Since the (pointwise) product
  $\prod_{j} v_{j}$ exists, we can choose a subsequence $(u_{i})$ of
  $(v_{j})$ satisfying
  $$ |1 - u_{i}(y^{-1},x)| \leq 1/i^2 \ \ \text{for all}\  x \in H_{i}, y
  \in F_{i}, i \in \Naturali .$$
  Let $x \in G$ and choose $N \in \Naturali$ such that $x \in H_{N}$.
  Then we get
  $$  \sum_{i} \frac{1}{\# F_{i}} \sum_{y \in F_{i}}|1 -
  u_{i}(y^{-1},x) | $$
  $$ \leq \sum_{i<N}2 \ + \ \sum_{i\geq N} \frac{1}{\# F_{i}} \sum_{y \in
  F_{i}} 1/i^2$$
  $$ = 2(N-1) \ + \ \sum_{i\geq N} 1/i^2 < \infty.$$
 This shows that $(F_{i})$ satisfies $(*)$ in Theorem 3.3,
 from which the result then clearly follows.
 \end{proof}

 \begin{corollary} Let $G$ be countable and amenable. Then there
 always exist  some  product sequence $(u_{i})$ in $Z^2(G,\Toro)$
 satisfying $u_{i} \neq 1$ for all $i$
 and   some sequence $\phi =
 (\phi_{i})$ of unit vectors in  $\ell^2(G)$ such that $\otimes_{i}
\lambda_{u_{i}}$
 exists on $\mathcal{H}^{\phi} = \otimes_{i}^{\phi}\ell^2(G)$. If
 $H^2(G,\Toro)$ is non-trivial and $1 \neq [u] \in H^2(G,\Toro)$,
 then the sequence $(u_{i})$ above may chosen so that $u = \prod_{i} u_{i}.$
 \end{corollary}

 \begin{proof}  We call a product sequence $(u_{i})$ in $Z^2(G,\Toro)$
 1-free if $u_{i} \neq 1$ for all $i$. It is easy to see that
 1-free product sequences
 do exist in $B^2(G,\Toro)$.
 As 1-freeness is clearly preserved when passing to subsequences,
 the first assertion  follows from the previous corollary.
 The  1-free product sequence $(u_{i})$  is then in $B^2(G, \Toro)$.
 Therefore (by closedness) $\prod_{i} u_{i} \in B^2(G, \Toro)$,
 so we may write it as ${\rm d} \rho$ for some normalized
 $\rho : G \rightarrow \Toro$.
 Assume now $H^2(G,\Toro)$ is non-trivial
 and $1 \neq [u] \in H^2(G,\Toro)$.
  Set $v_{1}= \overline{d\rho} \, u$ and $v_{i}= u_{i-1},i>1$.
 Then $(v_{i})$ is a 1-free product sequence satisfying
 $u = \prod_{i} v_{i}.$
 Further we can define a  new sequence $\psi = (\psi_{i})$ of unit vectors
 in $\ell^2(G)$, by setting $\psi_{1}= \delta_{e}$ and
 $\psi_{i}=\psi_{i-1}, i>1.$
 It is then obvious that $\otimes_{i} \lambda_{v_{i}}$ exists
 on $\mathcal{H}^{\psi}$, which proves the second assertion.
 \end{proof}

 \vspace{2ex}

 \noindent {\bf Remarks.} \\

 1) It follows from Theorem 3.1 iii) that representations obtained as
 the infinite tensor product of projective regular representations
 are never irreducible.

 \vspace{1ex}
 2) Let $G$ be countable and amenable, and let $(u_{i})$ and $(v_{i})$ be
two
sequences in
$Z^2(G,\Toro)$ satisfying $ v_{i}\sim_{\rho_{i}} u_{i}$ for every
$i$. Assume that $\otimes_{i} \lambda_{u_{i}}$
 exists on $\mathcal{H}^{\phi} = \otimes_{i}^{\phi}\ell^2(G)$ for some
 sequence $\phi = (\phi_{i})$ of unit vectors in  $\ell^2(G)$. As
 $\prod_{i} v_{i}$ does not necessarily exist, it may
  happen that $\otimes_{i} \lambda_{v_{i}}$ can not be formed at all
  (cf. Theorem 3.1). However, it is quite clear that
$\rho_1\lambda_{v_1} \otimes \rho_2\lambda_{v_2} \otimes \cdots$
exists on $\otimes^{\psi_i} \ell^2(G)$, where $\psi_{i}$ is defined by
$\psi_{i}(x) = \overline{\rho_{i}(x^{-1})} \phi_i(x)$, and this may be
considered as a problem of gauge fixing. On the other hand, let us also
assume that $\otimes_{i} \lambda_{v_{i}}$
 exists on $\mathcal{H}^{\psi} = \otimes_{i}^{\psi}\ell^2(G)$ for some
 sequence $\psi = (\psi_{i})$ of unit vectors in  $\ell^2(G)$. Then we
 may conclude that $\otimes_{i} \lambda_{v_{i}}$ is, up to unitary
 equivalence, just a perturbation of $\otimes_{i} \lambda_{u_{i}}$.

 (To prove this, we first appeal to Theorem 3.1 and obtain that both
 $u = \prod_{i} u_{i}$ and  $v = \prod_{i} v_{i}$ exist. Using Lemma 2.1
 we may  then write $v = {\rm d} \rho \, u$ for some normalized
 $\rho:G \rightarrow \Toro$. Now, using that $ \lambda_{v} \simeq \rho
 \lambda_{u}$ and Theorem 3.1, we get
 $$ \otimes_{i} \lambda_{v_{i}} \simeq \lambda_{v} \otimes id \simeq
 \rho ( \lambda_{u} \otimes id) \simeq \rho   \otimes_{i}
 \lambda_{u_{i}},$$
 where $id$ denotes the identity representation of $G$ on any infinite
 separable Hilbert space.)

 \vspace{1ex}

 3) To produce examples of infinite tensor product of projective unitary
  representations of not necessarily amenable groups, one can proceed
  as follows. Let $G$ be any countable group possessing a non-trivial
amenable factor group $K$ (one can here for instance let $G$ be any
non-perfect,
non-amenable group, e.\ g.\ any non-abelian countable free group, since
the abelianized group $G/[G,G]$ is
  then non-trivial and abelian) and let $(v_{i})$ be a sequence in
$Z^2(K,\Toro)$ such that
  $\otimes_{i} \lambda_{v_{i}}$
 exists on $\otimes_{i}^{\phi}\ell^2(K)$. Using
 the canonical homomorphism $\pi: G \rightarrow K$, we may lift each
 $v_{i}$ to a $u_{i} \in Z^2(G,\Toro)$ in the obvious way.
 Set $U_{i}(x):= \lambda_{v_{i}}(\pi(x)), x \in G,$  for each
 $i$. It is then a simple matter to check that each $U_{i}$ is a projective
unitary
  representation of $G$ on $\ell^2(K)$ associated to $u_{i}$, and
  that $\otimes_{i} U_{i}$ exists on $\otimes_{i}^{\phi}\ell^2(K)$.

\vspace{3ex}

We now turn to another class of examples which is in spirit related
to the setting of the Stone-Mackey-von Neumann theorem, i.\ e.\ with
so-called CCR-representations of a locally compact abelian group and its
dual
(cf. \cite{Sl}).

\vspace{1ex}
Let $A$ and $B$ be two discrete groups and $\sigma: A \times B \rightarrow
\Toro$ be
 a  bilinear map. We call a triple
 $\{V,W,\mathcal{H} \} $ for a \emph{CCR-representation} of $\sigma$
 whenever $V$ and $W$ are unitary representations of respectively $A$
 and $B$ on $\mathcal{H}$ which satisfy the CCR-relation
 $$V(a)W(b) = \sigma(a,b) \ W(b)V(a)$$
 for all $a \in A, b \in B.$

We now set $G= A \times B$ and  define $u_{\sigma}:
G \times G \rightarrow \Toro$ by
$$ u_{\sigma}(\,(a_{1},b_{1}),(a_{2},b_{2})\,) =
\overline{\sigma(a_{2},b_{1})}.$$
It is an easy exercise to check
that $ u_{\sigma}$ is a 2-cocycle on $G$ (in fact a bicharacter, i. e.
a bilinear map on $G \times G$ into $\Toro$).
When both $A$ and $B$ are abelian, then $[u_{\sigma}] \neq 1$ in
$H^2(G,\Toro)$ whenever $\sigma$ is non-trivial, as follows from \cite{Kl2}
since $u_{\sigma}$ is
then clearly non-symmetric. Note that there is an  1-1 correspondence
 between CCR-representations of $\sigma$ and projective unitary
representations of $G$ associated with $u_{\sigma}$ ( being
given by setting $U(a,b) = V(a) W(b)$ whenever
$\{V,W,\mathcal{H} \} $ is a CCR-representation of $\sigma$).

There is a canonical way to produce a CCR-representation of $\sigma $
on $\ell^2(B)$, to which we  may associate a projective unitary
representation
$U_{\sigma}$ of $G$ on $\ell^2(B)$ associated with  $ u_{\sigma}.$ We recall
this construction (and remark that a similar construction can be done on
$\ell^2(A)$ in
an analogous way):

For each
$a \in A, b\in B$ we set $\sigma_{a}(b) = \sigma(a,b)$, so  the map $(a
\mapsto
\sigma_a)
$ belongs to $ {\rm Hom}(A,\hat{B})$ where $\hat{B}:= {\rm Hom}(B,\Toro)$.
 Let then $V_{\sigma}(a)$ denote the multiplication operator by the
 function $\sigma_{a}$ on $\ell^2(B)$ and $\lambda_{B}$ be the left
 regular representation of $B$ on $\ell^2(B)$. By computation we have
 $$ V_{\sigma}(a)\lambda_{B}(b) = \sigma(a,b) \
 \lambda_{B}(b)V_{\sigma}(a)$$
 for all $a \in A, b \in B.$ Hence, the triple
$\{V_{\sigma},\lambda_{B},\ell^2(B) \} $ is
 a CCR-representation of $\sigma$ and we can put $U_{\sigma}(a,b) :=
 V_{\sigma}(a)\lambda_{B}(b)$ for all $(a,b) \in G$.

   \vspace{2ex}
   Assume now that $(\sigma_{i})$ is a sequence of bilinear maps from
  $ A \times B $ into $ \Toro.$ The question
whether is it possible to form  $\otimes_{i}U_{\sigma_{i}}$  on
  $\otimes_{i}^{\phi} \ell^2(B)$  for some sequence $\phi = (\phi_{i})$
  of unit vectors in $\ell^2(B)$
 is then  clearly equivalent  to whether it is possible to
  form the infinite tensor product of the CCR-representations
  associated with  the $\sigma_{i}$'s.  In the case of a positive
  answer, the product $\prod_{i}u_{\sigma_{i}}$ will  exist (as a
consequence of
  Theorem 3.1), so $\prod_{i}\sigma_{i}$ will then exist too and the
infinite tensor product of the CCR-representations
  associated with  the $\sigma_{i}$'s will be a CCR-representation of this
  product map.

  \vspace{1ex}
  Quite similarly to Theorem 3.2 and Theorem 3.3 we have:

  \newpage
  \begin{theorem}Let $(\sigma_{i})$ be a sequence of bilinear maps from
  $G = A \times B $ into $ \Toro.$ Set $U_{i}:= U_{\sigma_{i}}.$

  \vspace{1ex}\noindent
  i) Assume that $\otimes_{i}U_{i}$ exists on $\otimes_{i}^{\phi} \ell^2(B)$
  for some sequence $\phi = (\phi_{i})$ of unit vectors in $\ell^2(B).$
  Then $B$ is countable and amenable.

   \vspace{1ex}\noindent
  ii) Assume that $B$ is countable and amenable, and that $(F_{i})$ be a
$\sigma$F-sequence for
  $B$ satisfying $$\sum_{i} \frac{1}{\#(F_{i})} \ \sum_{b \in F_{i}} |1 -
  \sigma_{i}(a,b) | < \infty $$
  for every $a \in A.$ Set $\phi = (\phi_{i})$ where
$\phi_{i}:=\chi_{F_{i}}/{\#
 (F_{i}})^{1/2}.$

 \vspace{1ex}\noindent
 Then $\otimes_{i}U_{i}$  exists  on $\otimes_{i}^{\phi} \ell^2(B).$

  \end{theorem}

  \begin{proof} i) Since $U_{i}(e,b)= \lambda_{B}(b)$, this follows
  from \cite{ABC} (or Theorem 3.2).

  ii)  Let $B$ be countable and amenable, and $(F_{i})$ be as in ii). Since
  $(F_{i})$ is a  $\sigma$F-sequence for $B$ it follows from
  \cite{C2} (or Theorem 3.3) that $\otimes_{i} U_{i}(e,b) =
  \otimes_{i} \lambda_{B}(b)$ exists on  $\otimes_{i}^{\phi}
  \ell^2(B)$ for every $b \in B.$ The existence of $\otimes_{i}U_{i}$
  on $\otimes_{i}^{\phi} \ell^2(B)$ reduces then to whether
  $\otimes_{i}V_{\sigma_{i}}$ exists on
  $\otimes_{i}^{\phi} \ell^2(B)$, i.\ e.\ whether
  $$\sum_{i}|1-(V_{\sigma_{i}}(a)\phi_{i}, \phi_{i})| =
      \sum_{i}|1-((\sigma_{i})_{a} \,\phi_{i}, \phi_{i})| < \infty $$
  holds for every $a\in A$.
  As we have
  $$|1-((\sigma_{i})_{a} \,\phi_{i}, \phi_{i})|=  \frac{1}{\#(F_{i})} \
|\sum_{b \in F_{i}} (1 - \sigma_{i}(a,b))| \leq \frac{1}{\#(F_{i})} \
\sum_{b \in F_{i}} |1 - \sigma_{i}(a,b)| $$
  for every $a \in A$, this follows from the assumption on $(F_{i})$.
  \end{proof}

  We leave to the reader to deduce from this theorem the analogous versions
of
  Corollary 3.4 and Corollary 3.5 in this setting.

\section{The case of free abelian groups}

The purpose of this section is to examplify the results of the
previous section in the concrete case where G is a finitely generated
free abelian group.

\vspace{2ex}
\noindent We let $N \in \Naturali$ and set $G=\Relativi^N$.

\vspace{1ex}\noindent When $x=(x_1,\ldots,x_N) \in G$, we set
$|x|_1=\sum_{j=1}^N |x_j| \quad .$

\vspace{1ex}\noindent When $m \in \Naturali$, we define $K_m \subset G$ by
$$K_m=\{x \in G \ | \ 0 \leq x_{i} \leq m, i = 1 \ldots N \, \} \
(=\{0,1,\ldots,m\}^N).$$
To each $N \times N$ real matrix $A$, one may associate
$u_A \in Z^2(G,\Toro)$ by
$$u_A(x,y)=e^{ix\cdot(Ay)}.$$
In fact, every element in $H^2(G,\Toro)$ may be written as
$[u_A]$ for some skew-symmetric $A$ (see \cite{Bac,BB}).
Without loss of generality, we can assume that
$A \in M_N((-\pi,\pi])$,
i.e. all of $A$'s coefficients belong to $(-\pi,\pi]$. We set
$$|A|_\infty=\max\{|a_{ij}|, 1 \leq i,j \leq N \}.$$

\vspace{1ex}
We first record a technical lemma.

\bigskip
\begin{lemma}
 Let $A \in M_N((-\pi,\pi])$, $x,y \in G$ and $m \in \Naturali$. Then
\begin{itemize}
\item[(1)]
$|1-u_A(x,y)| \leq |A|_\infty |x|_1 |y|_1$
\item[(2)]
$\sum_{x\in K_m}|x|_1=\frac{Nm(m+1)^N}{2}$
\item[(3)]
$1-\frac{\#((x+K_m)\cap K_m)}{\#K_m} \leq \frac{|x|_1}{m+1} \ .$
\end{itemize}
\end{lemma}
\medskip

\begin{proof}

1) follows from $|1-e^{ix\cdot (Ay)}| \leq |x\cdot (Ay)|
\leq |A|_\infty |x|_1 |y|_1$.\par
\vspace{1ex}
2) $\sum_{x\in K_m}|x|_1
=\sum_{j=1}^N \sum_{x\in K_m}|x_j|
=N(m+1)^{N-1}(\sum_{k=0}^m k)=\frac{Nm(m+1)^N}{2}$.\par
\vspace{1ex}
3) $1-\frac{\#((x+K_m)\cap K_m)}{\#K_m}
=\frac{\#(K_m \backslash (x+K_m))}{\# K_m}
\leq \frac{(m+1)^{N-1}}{(m+1)^N}|x|_1
= \frac{|x|_1}{m+1}.$
\end{proof}

\medskip
\begin{proposition}
 Let $(A_i)$ be a sequence in $M_N((-\pi,\pi])$ and $(m_i)$ be a
sequence in $\Naturali$. For each $i \in \Naturali$, we set
$$F_i = K_{m_i} \subset G, $$
$$\phi_i=\frac{1}{(\# F_i)^{1/2}}\chi_{F_i} \in \ell^2(G), $$
$$u_i=u_{A_i} \in Z^2(G,\Toro).$$
Then we have:\par
\begin{itemize}
\item[(1)]
$(F_i)$ is a $F$-sequence for $G$ if and only if
$m_i \to +\infty$.
\item[(2)]
$(F_i)$ is a $\sigma F$-sequence for $G$ if and only if
$\sum_{i=1}^{\infty}\frac{1}{m_i} < \infty$.
\item[(3)]
$\prod_i u_i$ exists $\Leftrightarrow$ $\sum_i |A_i|_\infty < \infty$.
\item[(4)]
The projective unitary representation $\otimes_i \lambda_{u_i}$ of $G$
exists on $\otimes_i^{\phi_i}\ell^2(G)$
whenever $$\sum_{i=1}^{\infty}\frac{1}{m_i} < \infty \
\mbox{and} \ \sum_{i=1}^{\infty}m_i |A_i|_\infty< \infty$$
(and  $\prod_i u_i$ is then the associated 2-cocycle).
\end{itemize}
\end{proposition}
\medskip

\begin{proof}
The nontrivial parts of (1) and (2) are consequences of Lemma 4.1,
 part (3).
Assertion (3) relies on the inequality
$2|\theta|/\pi \leq |1-e^{i\theta}| \leq |\theta|$
which holds when $|\theta| \leq \pi$.
Concerning (4) let $x,y \in G.$ Then we have
\begin{align*}
\frac{1}{\# F_i}\sum_{y \in F_i} |1-u_i(-y,x)|
& \leq \frac{1}{(m_i+1)^N}(\sum_{y \in F_i}|A_i|_\infty |x|_1 |y|_1)
\quad \mbox{(by Lemma 4.1, (1))}
\\
& = \frac{|x|_1 |A_i|_\infty}{(m_i+1)^N} \sum_{y \in F_i} |y|_1 \\
& = \frac{|x|_1 |A_i|_\infty}{(m_i+1)^N} \frac{Nm_i(m_i+1)^N}{2}
\quad \mbox{(by Lemma 4.1, (2))}
\\
& = \frac{N|x|_1}{2}m_i|A_i|_\infty
\end{align*}
for every $i \in \Naturali$. Hence we have
$$\sum_i \frac{1}{\# F_i} \sum_{y \in F_i}|1-u_i(-y,x)|
\leq \frac{N|x|_1}{2}\sum_i m_i|A_i|_\infty .$$
Now if we assume that $\sum_{i=1}^{\infty}\frac{1}{m_i} < \infty$
and $\sum_{i=1}^{\infty}m_i |A_i|_\infty< \infty$,
then $\{F_i\}$ is a $\sigma F$-sequence for $G$ (by (2)) and
$\sum_i \frac{1}{\# F_i} \sum_{y \in F_i} |1-u_i(-y,x)| < \infty$
for all $x \in G$,
and the conclusion follows from Theorem 3.3.
\end{proof}

\medskip
\noindent {\bf Example.}
Let $A \in M_N((-\pi,\pi]).$ Set $A_i = 2^{-i}A$ and
$u_i=u_{A_i}$ ($i \in \Naturali$).
Then clearly $u_A=\prod_i u_{i}$.
Further, if we let $m_i=i^2$, then $\sum_i 1/m_i < \infty$ and
$\sum_i m_i|A_i|_\infty=|A|_\infty \sum_i i^2/2^i < \infty$
so (4) in the above proposition applies.
Theorem 3.1 then gives
$$\lambda_{u_A} \otimes id \cong \otimes_i \lambda_{u_i},$$
thus producing an infinite tensor product decomposition
of the amplification of $\lambda_{u_A}$.
It is well known that the C*-algebra $C^{*}(\lambda_{u_A})$ generated by
$\lambda_{u_A}$ on
$\ell^2(G)$ is a so-called non-commutative $N$-torus. Using this
decomposition result, we  can clearly obtain  a faithful representation of
$C^{*}(\lambda_{u_A})$ onto the
 $C^{*}$-algebra  generated by $\otimes_i \lambda_{u_i}$ on
 $\otimes_{i}^{\phi} \ell^2(G)$ for some suitably chosen sequence $\phi$
 of unit vectors in $\ell^2(G)$.

\vspace{2ex}
We shall now exhibit projective unitary representations arising from
CCR-representations of bilinear maps on some direct product
decomposition of $G$.

\vspace{1ex}
We assume from now on that $N\geq 2$ and write $G=\Relativi^N \simeq
\Relativi^P \times \Relativi^Q$ where $1 \leq P,Q < N$ and $ P+Q=N$.

To each $P \times Q$ matrix $D$ with coefficients in $(-\pi,\pi]$,
we associate a bilinear map
$\sigma_D: \Relativi^P \times \Relativi^Q \to \Toro$ by
$$\sigma_D(a,b)=e^{ia\cdot(Db)} .$$
Using the construction described at the end of the previous
section, we then obtain a CCR-representation of $\sigma_D$
on $\ell^2(\Relativi^Q)$, or, equivalently, a projective unitary
representation $U_D$ of $G=\Relativi^N$ with associated 2-cocycle $u^D$.
This cocycle is easy to describe:
a simple computation gives
$$u^D(x,y)=e^{ix\cdot(\tilde{D}y)} \quad (x,y \in G)$$
where $\tilde{D}$ is the $N \times N$ matrix given by
$$\tilde{D}=
\begin{pmatrix}
0 & 0 \\
-D^t & 0 \\
\end{pmatrix} .$$
Notice that $u^D=u_{\tilde{D}}$ and $[u^D]$ is non-trivial whenever $D
\neq 0$.

\medskip
 \begin{proposition}
 Let $(D_i)$ be a sequence of $P \times Q$ matrices
with coefficients in $(-\pi,\pi]$, and let $(U_i)=(U_{D_i})$
be the associated sequence of projective unitary representations
of $G$ on $\ell^2(\Relativi^Q)$.
Let $(n_i)$ be a sequence in $\Naturali$.

\vspace{1ex}\noindent Set $H_i=\{b \in \Relativi^Q \ | \ 0 \leq b_{i} \leq
n_i,
i = 1 \ldots Q \}$
and $\psi_i=1/(\# H_i)^{1/2} \chi_{H_i} \ (i \in \Naturali)$.  \\
\vspace{1ex}\noindent Then $\otimes_i U_i$ exists on
$\otimes_i^{\psi_i}\ell^2(\Relativi^Q)$
whenever $\sum_i 1/n_i < \infty$ and $\sum_i n_i |D_i|_\infty < \infty.$
\end{proposition}

\begin{proof}
This follows from Theorem 3.6. As the details are quite similar to the
proof of the previous proposition, we leave these to the reader.
\end{proof}

\medskip
\noindent {\bf Example.}
We take $P=Q=1$ so that $G=\Relativi \times \Relativi=\Relativi^2$,
and let $(D_j)=(\theta_j)$ be a sequence in $(-\pi,\pi]$.
This gives rise to the sequence $(U_j)$ of representations of
$\Relativi^2$ on $\ell^2(\Relativi)$ with associated 2-cocycles
$$u_j(x,y)=e^{-i\theta_j x_1 y_2} \ (x,y \in \Relativi^2) .$$
By Proposition 4.3 we can then form the infinite tensor representation
$\otimes_j U_j$ whenever we can choose a sequence $(n_j)$ in $\Naturali$
such that $\sum_j 1/n_j < \infty$ and $\sum_j n_j |\theta_j| <\infty $
(e.g. $n_j=j^2$ will do  if $(j^4|\theta_j|)$ is bounded).

By a more careful analysis of this example
involving
the familiar
Dirichlet sums,
one can deduce that $\otimes_j U_j$ will exist whenever
we can choose $(n_j)$ such that
$$\sum_j \frac{1}{n_j}<\infty \ \text{and} \
\sum_j |1-\frac{1}{2n_j +1}
\frac{\sin((2n_j +1)\theta_j/2)}{\sin(\theta_j/2)} \, | < \infty .$$
Assuming that $\sum_j |\theta_j|<\infty$ (so $\prod_j u_j$ exists), it
would be interesting to know whether such a choice  of $(n_{j})$ can
always be made.

\section{Infinite products of actions}

For each $i \in \Naturali$
let $\H_i $ be a Hilbert space, $\phi_i \in \H_i$  be a  unit vector,
$\M_i \subset \B(\H_i)$ be a von Neumann algebra
and $\alpha_i: G \to  {\rm Aut}(\M_i)$  be an action of $G$ on $\M_{i}$.
We denote by $I_i$  the identity operator on $\H_i$.
We then form the  $*$-algebra $\odot_i\M_i$ (resp. von Neumann
algebra  $\otimes_i(\M_i,\phi_i)$) acting on $\otimes_{i}^{(\phi_{i})}
\H_{i}$ generated by operators of the form $ \otimes_{i} T_{i}$ where
$T_{i} \in \M_{i}$ and $T_{i} = I_{i}$ for all but finitely many $i$'s.
At the $*$-algebraic level we define  an action
$\odot_{i}\alpha_i$ of $G$ on $\odot_i \M_i$
such that for every finite $J \subset \Naturali$ we have
$$\odot_{i}
 \alpha_i
((\otimes_{i \in J} T_i)\otimes(\otimes_{i \notin J} I_i))
= (\otimes_{i \in J} \alpha_i(T_i))\otimes(\otimes_{i \notin J}
I_i).$$
One natural question is whether $\odot_{i}
\alpha_i$ may be extended to an action of $G$ on the von Neumann
algebra
$\otimes_i(\M_i,\phi_i)$. As we shall see, the answer may be negative in
some situations, regardless of the choice of unit vectors $\phi_{i}$.

\vspace{1ex}
We  retrict ourselves to the case where each
$\alpha_{i}$ is unitarily implemented, i.\ e.\ we assume that for every $i$
and
$g$ there exists a unitary $ U_i(g)$ on $\H_i$
such that $\alpha_{i,g}={\rm Ad}\,(U_i(g)).$
This assumption is
automatically satisfied
for many classes of von Neumann algebras (see \cite{StZs}, \S 8).
Note that if $U_{i}(g) \in \M_{i}$ for all $g \in G$ and $M_{i}$ is a
factor, especially if $\M_{i} = \B(\H_{i}),$ then $ g \to U_{i}(g)$ is
a projective unitary representation of $G$ on $\H_{i}.$

We consider the following condition:
$$(*) \ \ \sum_i (1- |(U_i(g)\phi_i,\phi_i)|) < \infty \quad\quad
\text{for all}\  g \in G.$$

\begin{proposition} Condition $(*)$ is equivalent to the following
condition:
$$ (**) \ \exists \ \rho_i: G \to \Toro, \rho_i(e)=1, \
\mbox{such that} \ \otimes_i \rho_i U_i \ \mbox{exists on} \
\otimes_i^{\phi_i} \H_i.$$
When $(*)$ holds, then $\odot_{i} \alpha_i$
extends to a
unitarily implemented action $\alpha$ on $\otimes_i(\M_i,\phi_i),$
which is inner whenever $U_i(g) \in \M_i$ for every $i$
and $g \in G.$
\end{proposition}

\begin{proof} The first assertion follows from Proposition 2.3, using
\cite[\S 1.2]{Gui}. When $(*)$ holds, then $\alpha_g = {\rm Ad}\,(U(g))$
where $U(g)=\otimes_i\rho_i(g)U_i(g)$ is well defined on
$\otimes_i^{\phi_i}\H_i.$
 Clearly $U(g) \in \otimes_{i}(\M_i,\phi_{i})$ whenever $U_i(g) \in
\M_i$ for every $i$
and $g \in G,$
and $\alpha_g$ is then inner for every $g \in G$.
\end{proof}

We now treat the case where every $\M_{i}$ is a type $I$ factor. We use the
well known fact that every automorphism of a
type $I$ factor is inner and also that $\otimes_i(\B(\H_i),\phi_i)
\ = \B(\otimes_i^{\phi_i}\H_i)$ ( \cite[Proposition 1.6]{Gui}).

\begin{theorem} Assume that $\M_{i} = \B(H_{i})$ for all $i.$ Then
$\odot_{i} \alpha_i$ extends
(uniquely) to an action $\alpha =\otimes \alpha_i$ on
$\otimes_i(\B(\H_i),\phi_i)$
if and only if condition $(*)$ holds.
\end{theorem}

\begin{proof}
Assume that an extension $\alpha$ of $\odot_{i} \alpha_i$
exists on $\M^{\phi} =\otimes_i(\M_i,\phi_i).$ Using the facts recalled
above,
we have $\alpha_{g} = {\rm Ad}\,(U(g))$
for some  $U(g) \in {\cal U}(\otimes_i^{\phi_i}\H_i)$ for every $g \in G$.
\\
Let $J$ be a non-empty finite subset of $\Naturali$.\\
We identify
$\M^{\phi}$ with $(\otimes_{i \in J} \M_i) \otimes \, {}_{J}\M$
where ${}_{J}\M :=\otimes_{i \notin J}(\M_i,\phi_i)$,
and  consider ${}_{J}\M $ as a von Neumann subalgebra of
$\M^{\phi}$ in the obvious way. It is easy to see that $\alpha$
restricts to an action ${}_J\alpha$ of $G$ on ${}_{J}\M $ such that
$\alpha = (\otimes_{i \in J}\alpha_i)\otimes {}_J\alpha$.
Since${}_{J}\M$ is a also type $I$ factor, we can write
${}_J\alpha_{g}={\rm Ad}\,({}_JU(g))$ for some
${}_JU(g) \in {\cal U}(\otimes_{i \notin J} ^{\phi_i}\H_i)$ for each $g \in
G.$

Set now $U_J(g)=\otimes_{i \in J}U_i(g)$ for each $g \in G$.
Then
$\alpha_g={\rm Ad}\,(U_J(g)\otimes {}_JU(g))$.
Therefore, for each $g \in G$, there exists some $z_{J}(g) \in \Toro$ such
that
 $U(g)=z_{J}(g)U_J(g)\otimes {}_JU(g)$.

 Let $g \in G.$ Since $U(g) \neq 0$ we can pick two elementary decomposable
vectors
 $\otimes \psi_{i}$ and $\otimes \xi_{i}$ in
 $\otimes_i^{\phi_i}\H_i$ (which
do not depend on $J$) satisfying
$$0 \neq c(g):=|(U(g) \otimes \psi_{i} , \otimes \xi_{i})|
=\prod_{i \in J}|(U_i(g)\psi_i , \xi_i)|
\ |({}_JU(g)\otimes_{i \notin J}\psi_i,\otimes_{i \notin J}\xi_i)|$$
Since $|({}_JU(g)\otimes_{i \notin J}\psi_i,\otimes_{i \notin
J}\xi_i)| \leq 1$ we get
$$ 0 < c(g) \leq \prod_{i \in J}|(U_i(g)\psi_i , \xi_i)|.$$
As this holds for every $J$, one easily deduces  that $\prod_{i \in
\Naturali}|(U_i(g)\psi_i , \xi_i)|$ converges
to a non-zero number. Since $\psi_{i} =\xi_{i} = \phi_{i}$ for all but
finitely many $i$'s, this implies that $(*)$ holds. Hence, we have
shown the only if part of the assertion. The converse part follows from
Proposition 5.1.
\end{proof}

The proof of the above result is reminiscent of the proof of
a lemma  in \cite{Sto} (see also \cite{El}). In the same line of
ideas, we have the following  result, which is
related to \cite[Lemme 1.3.8]{Co}.

\begin{theorem} Assume that all $\M_i$'s are factors and that
$\odot_{i} \alpha_i$ extends to an action $\alpha$ on $\M^\phi
=\otimes_i(\M_i,\phi_i).$
Then $\alpha$ is inner if and only if there exists for each $g \in G$
and each $i$ a  unitary $v_i(g) \in \M_i$  implementing $\alpha_{i,g}$
such that the following condition holds:
 $$ (1) \ \ \sum_i (1- |(v_i(g)\phi_i,\phi_i)|) < \infty \quad\quad
\text{for all}\  g \in G.$$ On the other hand, $\alpha$ is outer if and
only if, for each $g \in G, g\neq e,$ at least one of the $\alpha_{i,g}$
is outer or there exists for each $i$ a  unitary $v_i(g) \in \M_i$
implementing
$\alpha_{i,g}$  such that
 $$ (2) \ \ \sum_i (1- |(v_i(g)\phi_i,\phi_i)|) = \infty.$$
\end{theorem}

\begin{proof} Assume first that  $\alpha$ is inner. So we
have $\alpha_{g} = {\rm Ad}\,(U(g))$
for some  unitary $U(g) \in \M^\phi$ for every $g \in G$. Recall
from \cite{Gui} that $\M^{\phi}$ is a factor.
Using \cite[Corollary 1.14]{Ka}, it follows easily that each $\alpha_{i}$ is
inner. Hence, there exists for each $g \in G$
and each $i$ a  unitary $v_i(g) \in \M_i$  implementing $\alpha_{i,g}.$

Let $J$ be a non-empty finite subset of $\Naturali$. As in the
previous proof, we identify
$\M^{\phi}$ with $(\otimes_{i \in J} \M_i) \otimes \, {}_{J}\M$
where ${}_{J}\M :=\otimes_{i \notin J}(\M_i,\phi_i)$,
We set $V_J(g)=\otimes_{i \in J}v_i(g)$ for each $g \in G$  and
$W_{J}(g)=(\, V_J(g) \otimes (\otimes_{i \notin J} I_i)\,)^*U(g).$
Then, using that we may write
$\alpha = (\otimes_{i \in J}\alpha_i)\otimes {}_J\alpha$, we get
 $$W_{J}(g) \in  ( \otimes_{i}(\M_i,\phi_i)) \cap ((\otimes_{i \in
J}\M_i) \otimes (\otimes_{i \notin J}\Complessi I_i))'.$$
Using that all $\M_i$ are factors, it is  a simple exercise to
deduce  that $W_{J}(g) \in
(\otimes_{i \in J}\Complessi I_i) \,\otimes \, (\otimes_{i \notin
J}(\M_i,\phi_i))$.
We may therefore write $W_{J}(g) = (\otimes_{i \in J} I_i) \otimes
{}_JV(g)$ for some unitary ${}_JV(g) \in \otimes_{i \notin
J}(\M_i,\phi_i)$. This gives $U(g) = V_J(g) \otimes {}_JV(g)$
and we can clearly proceed further in the same way  as in the
previous proof to show that  $(1)$ holds, thereby proving the only
if part of the first assertion. The converse part of this assertion follows
from Proposition 5.1.
The second assertion follows from a similar argument.
\end{proof}

The following corollary may be seen as generalization of
\cite[Theorem 6.7]{EK}.

\begin{corollary} Assume  for each $i \in \Naturali$ that  $\beta_{i}$ is
an action of $G$
on some von Neumann algebra $\N_{i}$ and that there exists  a
 normal $\beta_{i}$-invariant state $\tau_{i}$ on $\N_{i}.$ Denote the
GNS-triple of $\tau_{i}$ by $(\pi_{i}, \H_{i}, \xi_{i})$ and set
$\M_{i} = \pi_{i}(\N_{i}).$ Let $\alpha_{i}$ be the action of $G$
on $\M_{i}$ induced by $\beta_{i}.$ Then $\odot_{i} \alpha_{i}$
extends to an action $\alpha$ of $G$ on $\otimes_i(\M_i,\xi_i).$

Assume further that all $\N_{i}$'s are factors and all $\pi_{i}$'s are
faithful. Then $\alpha$ is inner
if and only if there exists for each $g \in G$
and each $i$ a  unitary $v_i(g) \in \N_i$  implementing $\beta_{i,g}$
such that the following condition holds:
 $$ (1)\ \ \sum_i (1- |\tau_{i}(v_i(g))|) < \infty \quad\quad
\text{for all}\  g \in G.$$ On the other hand, $\alpha$ is outer if and
only if, for each $g \in G, g \neq e,$ at least one of the $\beta_{i,g}$ is
outer
or there exists  each $i$ a  unitary $v_i(g) \in \N_i$  implementing
$\beta_{i,g}$  such that
 $$  (2) \ \ \sum_i (1- |\tau_{i}(v_i(g))|) = \infty.$$

\end{corollary}

\begin{proof} We first recall
that there exists for each $i$  a unitary representation $V_{i} :
G \to \B(\H_{i})$ such that
$$ \pi_{i}(\beta_{i,g}(x)) = V_{i}(g) \pi_{i}(x) V_{i}(g)^{*} \
\text{and} \  V_{i}(g) \pi_{i}(x) \xi_{i} = \pi_{i}(\beta_{i,g}(x)) \xi_{i}
$$
for all $g \in G, x \in \N_{i}$ (see \cite{Di}). The induced action
$\alpha_{i}$ on $\M_{i}$ is then defined by
$\alpha_{i,g}(\pi_{i}(x)) = \pi_{i}(\beta_{i,g}(x)).$
As $ V_{i}(g)  \xi_{i} =  \xi_{i}$ for all $g \in G$, the first
assertion follows obviously from Proposition 3.1. The second
assertion is then easily deduced from Theorem 5.3.
\end{proof}

\medskip
\noindent {\bf Example.} Let $u_{i}$ be  a sequence in $Z^2(G,\Toro).$
Set $\N_{i} = \lambda_{u_{i}}(G)^{''} \subset \B(\ell^2(G))$ and let
$\beta_{i,g}$ be the inner automorphism of $\N_{i}$ implemented by
$\lambda_{u_{i}}(g)$ for all $g \in G, i \in \Naturali .$ Let
$\tau_{i}$ denote the canonical normal faithful tracial state of
$\N_{i}$ (determined
by $\tau_{i}(\lambda_{u_{i}}(g)) = 1$ if $g = e$ and $0$ otherwise),
which is trivially $\beta_{i}$-invariant.
If $\xi$ denote the normalized delta-function at $e$, then
$\tau_{i}=\omega_{\xi}|_{\N_{i}}.$ So we may identify the GNS-triple
of $\tau_{i}$ with $(id_{i},\ell^2(G),\xi_{i})$, where $id_{i}$ denotes
the identity representation of $\N_{i}$ and $\xi_{i}=\xi,$ i.\ e.\
we may take $\M_{i}=\N_{i}$ and $\alpha_{i}= \beta_{i}$ in the notation of
Corollary 5.4. Hence, $\odot \alpha_{i}= \odot \beta_{i}$ extends to
an action $\alpha$ on $\otimes_{i} (\lambda_{u_{i}}(G)^{''}, \xi_{i}).$

Further, if all $\lambda_{u_{i}}(G)^{''}$ are factors, then $\alpha$
is outer, as $$\sum_i (1- |\tau_{i}(\lambda_{u_{i}}(g))|) = \sum_{i}
1 = \infty \ \ \text{for all} \ g \neq e.$$ A necessary and sufficient
condition
for a twisted group von Neumann algebra $\lambda_{u}(G)^{''}$ to be a factor
may be found in \cite{Kl1}.

If we replace each $\N_{i}$ with $\B(\ell^2(G))$ in this example, the
extended product action may be formed in many cases under the
assumption that $G$ is
countable and amenable, as follows from Teorem 3.3  and  Proposition 5.1.
This requires a suitable choice of unit vectors $\phi_{i}$ in
$\ell^2(G).$ This product action restricts then to an action on
$\otimes_{i} (\lambda_{u_{i}}(G)^{''}, \phi_{i})$ which is inner, in
contrast to the factor case above.
When $G$ is either uncountable or non-amenable, we have the following:

\begin{theorem}
 Let $u_i$ be a sequence in $Z^2(G,\Toro)$ and
$\alpha_i = {\rm Ad}\,\lambda_{u_i}$ be the associated
sequence of actions of $G$ on $\B(\ell^2(G)).$
If $G$ is either uncountable or non-amenable, then
$\odot_{i} \alpha_i$
does not extend to
an  action of $G$ on
$\otimes_i (\B(\ell^2(G)),\phi_i),$
regardless of the choice of vectors $\phi_i$.
\end{theorem}

\begin{proof}
According to  Proposition 5.1 and Theorem 5.2,
the existence of such an extension $\otimes_i (\B(\ell^2(G)),\phi_i)$ would
imply the
existence of $\otimes_i \rho_i \lambda_{u_i}$
on $\otimes_i^{\phi_i}\ell^2(G)$
for some choice of functions $\rho_i: G \to \Toro$ with $\rho_i(e)=1$.
It is straightforward to see that
this amounts to the existence of
$\otimes_i \lambda_{v_i}$ on  $\otimes_i^{\psi_i}\ell^2(G)$
for some $v_i \in Z^2(G,\Toro)$ with $v_i \sim u_i$ and some sequence
$\psi_{i}$ of unit vectors in $\ell^2(G).$
This is impossible if $G$ is either uncountable or non-amenable, as follows
from
Theorem 3.2.
\end{proof}

Another type of possible obstruction for extending a product action from
the $*$-algebraic level to the von Neumann algebra level is of
cohomological nature, as we now illustrate:

\begin{theorem}
Let $\alpha_i$ be a sequence of actions of  $G$ on $\B(\H_i)$
and write each $\alpha_i$ as ${\rm Ad}\,U_i(g)$
where $U_i$ is a projective representation of $G$ with
associated  2-cocycle $u_i.$
Assume that $[u_i] = [u]$ for every $i$ and $[u]\neq [1]$ in
$H^2(G,\Toro)$. \\
Then $\odot_{i} \alpha_i$ does not extend to an action of $G$ on
$\otimes_i(\B(\H_i),\phi_i),$ regardless of the choice of vectors
$\phi_i.$
\end{theorem}

\begin{proof}
Assume that such an extension exists $\otimes_i(\B(\H_i),\phi_i).$
Using Proposition 5.1 and Theorem 5.2, we deduce that
$\otimes \rho_i U_i$ exists on $\otimes_i^{\phi_i}\H_i$ for some choice
of functions $\rho_i: G \to \Toro$ with $\rho_i(e)=1$.
It follows then from Theorem 3.1 that  $\prod_i ({\rm d}\rho_i) u_i$
exists.
Hence
${\rm d}\rho_i u_i \to 1$ (in the pointwise topology). As each $u_{i}=({\rm
d}\rho'_i) u$  for some $\rho'_{i}$, we get that
 $u$ is a limit of 2-coboundaries. Since $B^2(G,\Toro)$ is
 closed, this means that $u$ is itself a
 coboundary, i.\  e.\ $[u]=1$, which gives  a contradiction.
\end{proof}

\medskip
\noindent {\bf Example.}
The simplest case where the above situation  occurs is when
 $G=\Relativi_2 \times \Relativi_2$. Indeed, let
$$V=
\begin{pmatrix}
0 & 1 \\
1 & 0 \\
\end{pmatrix}, \quad
W=\begin{pmatrix}
1 & 0 \\
0 & -1 \\
\end{pmatrix}. $$
A projective unitary representation of $G=\Relativi_2 \times
\Relativi_2$
on $\Complessi^2$
is then obtained by setting $U((a,b))=V^a \, W^b \ (a,b \in \Relativi_2)$.
Since $V^a \, W^b=\sigma(a,b)W^b \, U^a$ where $\sigma(a,b)=-1$ if $a=b=1$
and $1$ otherwise,
the associated cocycle $u$ is easily computed to be
$u((a_1,b_1),(a_2,b_2))=(-1)^{a_2b_1}.$
It is not difficult to check that $[u] \neq 1.$ Remark that $U$ is nothing
but the projective representation associated to the CCR
representation of $\sigma$ on $\Complessi^2 = \ell^2(\Relativi_2) $
determined by $V$ and $W.$

For each $i \in \Naturali$ consider the action $\alpha_{i}$ of $G$
on $M_2(\Complessi)$
given by $\alpha_{i,(a,b)}={\rm Ad}\,(U((a,b)))$.
Then, according to Theorem 5.6, the infinite tensor product
of the $\alpha_{i}$ 's does never make sense as an action on
 $\otimes_i (M_2(\Complessi),\phi_i)$.

On the other hand, the canonical tracial state of $M_2(\Complessi)$ is
trivially $\alpha_{i}$-invariant. Therefore we may use Corollary 5.4
to form the infinite tensor product action after passing to the
GNS-representation with respect to this tracial state for each $i.$
As another application of Corollary 5.4, the resulting
product action is easily seen to be outer.

\bigskip

\noindent{\bf Acknowledgements.}
Part of this work was made while R.\ C.\ was visiting the
Department of Mathematics at the University of Oslo
 during the academic year 1998-1999 supported by the EU TMR network
 ``Non Commutative Geometry''.
He is very grateful to the members of the Operator Algebra group
in Oslo for their friendly
hospitality
and for providing perfect working conditions.
Both authors would like to thank Ola Bratteli and Erling St\o rmer for
some helpful discussions. 

\bigskip
{\parindent=0pt Addresses of the authors:

\smallskip Erik B\'edos, Institute of Mathematics, University of
Oslo, \\
P.B. 1053 Blindern, 0316 Oslo, Norway. E-mail: bedos@math.uio.no. \\

\smallskip \noindent
Roberto Conti, Dipartimento di Matematica, Universit\`a di Roma  \\
``Tor Vergata'',
I-00133 Roma, Italy. E-mail: conti@mat.uniroma2.it.\par}


\vskip1cm

\bigskip


\begin{thebibliography}{40}


\bibitem{ABC} M. Aita, W. Bergmann, R. Conti: Amenable groups
and generalized Cuntz algebras,
{\it J. Funct. Anal.\/} {\bf 150}, 48--64 (1997).


\bibitem{Bac} N. B. Backhouse : Projective representations
of space groups, II: Factor systems, {\it Quart.\ J. Math.\ Oxford \/}
{\bf 21} (1970), 223-242.

\bibitem{BB} N. B. Backhouse, C. J. Bradley : Projective representations
of space groups, I: Translation groups, {\it Quart.\ J. Math.\
Oxford \/} {\bf 21} (1970), 203-222.


\bibitem{BC} W. Bergmann, R. Conti: Group representations and
generalized Cuntz algebras, in preparation.

\bibitem{Br} K. S. Brown: Cohomology of groups, GTM in Mathematics 87,
Springer-Verlag, New-York (1982).

\bibitem{Co} A. Connes: Une classification des facteurs de type III,
{\it Ann. Ec. Norm. Sup. \/} {\bf 6}, 133--252 (1973).

\bibitem{C2} R. Conti: Amenability and infinite tensor products,
in Operator Theory, Operator Algebras and Related Topics,
Proceedings of the 16th International Conference on
Operator Theory, Timisoara 1996, The Theta Foundation, Bucharest 1997.


\bibitem{Di} J. Dixmier: Les C*-algebres et leurs representations,
Gauthiers-Villars, Paris (1969).

\bibitem{El} G. A. Elliott:
Finite projections in tensor products of von Neumann algebras,
{\it Trans. Amer. Math. Soc.\/} {\bf 212}, 47--60 (1975).

\bibitem{EK} D. E. Evans, Y. Kawahigashi: Quantum symmetries on
operator algebras,  Oxford : Oxford University Press (1998).



\bibitem{Gui} A. Guichardet:
Produits tensoriels infinis et repr\'esentations des relations
d'anticommutation,
{\it Ann. Scient. \'Ec. Norm. Sup.\/} {\bf 3}, 1--52 (1966).

\bibitem{Gui2} A. Guichardet: Tensor product of C*-algebras, part  I
and II, {\it Lect. Notes Series No} {\bf 12,13}, Aarhus University (1969).

\bibitem{HiSh} T. Hirai, H. Shimomura:
Relations between unitary representations of diffeomorphism groups and
those of the infinite symmetric group or of related permutation groups,
{\it J. Math. Kyoto Univ.\/} {\bf 37}, 261--316 (1997).


\bibitem{Ka} R. R. Kallman: A generalization of free action,
{\it Duke Math. J. \/} {\bf 36}, 781--789 (1969).

\bibitem{Kl1} A. Kleppner: The structure of some induced representations,
{\it Duke Math. J.\/} {\bf 29}, 555--572 (1962).

\bibitem{Kl2} A. Kleppner: Multipliers on Abelian groups,
{\it Math. Ann.\/} {\bf 158}, 11--34 (1965).

\bibitem{vN} J. von Neumann:
On infinite direct products,
{\it Compos. Math.\/} {\bf 6}, 1--77 (1938).


\bibitem{Pat} A. Paterson:
Amenability, Math. Surveys and Monographs 29, Amer. Math. Soc.
(1988).

\bibitem{Pi} J. P. Pier:
Amenable locally compact groups,
John Wiley \& Sons Inc. (1984).


\bibitem{Sl} J. Slawny:
On factor representations and the C*-algebra of canonical
commutation relations,
{\it Comm. Math. Phys.\/} {\bf 24}, 151--170 (1972).

\bibitem{StZs} S. Str\u{a}til\u{a}, L. Zsid\'{o}:
Lectures on von Neumann algebras,
Abacus Press, Tunbridge Wells, Kent (1979).

\bibitem{Sto} E. St\o rmer:
On infinite tensor products of von Neumann algebras,
{\it Amer. J. Math.\/} {\bf 93}, 810--818 (1971).

\end{thebibliography}
\end{document}